\documentclass[a4paper,11pt]{amsart}
\usepackage[colorlinks, linkcolor=blue,anchorcolor=Periwinkle,
    citecolor=blue,urlcolor=Emerald]{hyperref}
\usepackage[all]{xy}
\SelectTips{cm}{}

\newcommand{\ie}{{\em i.e.}\ }

\newcommand{\ko}{\: , \;}

\newcommand{\ol}[1]{\overline{#1}}

\numberwithin{equation}{subsection}
\newtheorem{theorem}[subsection]{Theorem}
\newtheorem{classification-theorem}[subsection]{Classification Theorem}
\newtheorem{decomposition-theorem}[subsection]{Decomposition Theorem}
\newtheorem{proposition-definition}[subsection]{Proposition-Definition}
\newtheorem{periodicity-conjecture}[subsection]{Periodicity Conjecture}

\newtheorem{remark}[subsection]{Remark}

\newcommand{\reminder}[1]{}

\newcommand{\iso}{\xrightarrow{_\sim}}

\newcommand{\Hom}{\mathrm{Hom}}

\newcommand{\RHom}{\mathrm{RHom}}

\newcommand{\ten}{\otimes}
\newcommand{\lten}{\overset{\boldmath{L}}{\ten}}

\newcommand{\ca}{{\mathcal A}}
\newcommand{\cb}{{\mathcal B}}

\newcommand{\cd}{{\mathcal D}}

\newcommand{\cR}{{\mathcal R}}

\newcommand{\bp}{\mathbf{p}}

\renewcommand{\phi}{\varphi}

\renewcommand{\tilde}[1]{\widetilde{#1}}

\begin{document}

\date{June 3, 2018}

\title{Erratum to ``Deformed Calabi--Yau Completions''}
\author{Bernhard Keller}
\address{Universit\'e Paris Diderot -- Paris 7\\
    UFR de Math\'ematiques\\
   Institut de Math\'ematiques de Jussieu--PRG, UMR 7586 du CNRS \\
   Case 7012\\
    B\^{a}timent Chevaleret\\
    75205 Paris Cedex 13\\
    France
}
\email{bernhard.keller@imj-prg.fr}
\urladdr{https://webusers.imj-prg.fr/~bernhard.keller/}

\begin{abstract}
We correct an error and some inaccuracies 
that occurred in ``Deformed Calabi--Yau completions''. The most important point is that,
as pointed out by W.~K.~Yeung,
to show that the deformed Calabi-Yau completion has
the Calabi--Yau property, one needs to assume that the deformation parameter
comes from negative cyclic homology. Notice that this does hold in the case
of Ginzburg dg algebras.
\end{abstract}




\maketitle

\section{Calabi--Yau completions}
\label{s:CY-completions}

We refer to \cite{Keller11b} for unexplained notation and terminology.
Let $k$ be a commutative ring, $n$ an integer and $\ca$ a dg $k$-category. 
We may and will assume that $\ca$ is cofibrant
over $k$, \ie each morphism complex $\ca(X,Y)$ is cofibrant as a complex
of $k$-modules. Moreover, we assume that $\ca$ is homologically smooth
\cite{KontsevichSoibelman06}, 
\ie $\ca$ is perfect in the derived category $\cd(\ca^e)$ of $\ca$-bimodules. For a 
dg $\ca$-bimodule $M$, put
\[
M^\vee=\RHom_{\ca^e}(M,\ca^e).
\]
Recall \cite{Ginzburg06} that $\ca$ is {\em (bimodule) $n$-Calabi--Yau} if there is an isomorphism
\[
\phi: \Sigma^n \ca^\vee \iso \ca
\]
in the derived category of $\ca$-bimodules. 
The symmetry property originally imposed on $\phi$ is automatic as shown in
Appendix~C of \cite{VandenBergh15}.  Following 
\cite{ThanhofferVandenBergh12, KontsevichVlassopoulos13, Yeung18},
we define an {\em $n$-Calabi-Yau structure on $\ca$}
as the datum of a class $\eta$ in negative cyclic homology $HN_n(\ca)$ 
which is {\em non degenerate}, \ie whose
image under the canonical maps
\[
HN_n(\ca) \to HH_n(\ca,\ca) \iso \Hom_{\cd(\ca^e)}(\Sigma^n \ca^\wedge, \ca)
\]
is an isomorphism. Following \cite{VandenBergh15}, such a structure
is called {\em exact} if $\eta$ is an image under Connes' map
$B: HC_{n-1} (\ca) \to HN_n(\ca)$.

Let $\Theta_\ca$ be the inverse dualizing complex $\ca^\vee$ and
$\theta$ a cofibrant replacement of $\Sigma^{n-1}\Theta_\ca$.
The {\em $n$-Calabi--Yau completion $\Pi_n(\ca)$} was defined in \cite{Keller11b}
as the tensor category $T_\ca(\theta)$. The following theorem
is a more precise version of Theorem~4.8 of \cite{Keller11b}.
We include a proof since the statement is slightly
stronger and the new proof more transparent.

\begin{theorem} \label{thm:cy-completion}
The Calabi-Yau completion $\Pi_n(\ca)$ is homologically
smooth and carries a canonical exact $n$-Calabi-Yau structure.
\end{theorem}

\begin{remark} For finitely cellular dg categories $\ca$, 
W.~K.~Yeung gives two proofs of this theorem: one in section~3.3 of \cite{Yeung16} and a
more geometric one in section~2.3 of \cite{Yeung18}. The theorem generalizes his result 
to arbitrary homologically smooth dg categories $\ca$.
\end{remark}

\begin{remark} It is well-known that the Calabi-Yau completion,
being a generalization of the $2$-Calabi--Yau preprojective algebra \cite{CrawleyBoevey00}, should be viewed as a (shifted) non commutative 
cotangent bundle. This is very nicely explained in section~2.3 of \cite{Yeung18}.
Alternatively, one may view it as a (shifted) non commutative
total space of the canonical bundle.
This is made rigorous in section~3.5 of \cite{IkedaQiu18}.
\end{remark}

\begin{proof}[Proof of the Theorem] Put $\cb=\Pi_n(\ca)$. Notice that $\cb$ is
augmented over $\ca$ in the sense that we have canonical dg functors
$\ca\to\cb\to\ca$ whose composition is the identity. Thus, Hochschild and cyclic
homology of $\ca$ are canonically direct summands of those of $\cb$.
We call the supplementary summands the {\em reduced} Hochschild respectively
cyclic homology of $\cb$. We would like to compute them.
Since $\cb$ is a tensor category (though over the non commutative ground
category $\ca$), it suffices to adapt the results of section~3.1 of \cite{Loday98}. 
The analogue of the
small resolution $C^{\mathrm{sm}}(T(V))$ of remark 3.1.3 of \cite{Loday98}
is the exact bimodule sequence
\begin{equation} \label{eq:small-resolution}
\xymatrix{0 \ar[r] & \cb\ten_\ca\theta\ten_\ca\cb \ar[r]^-{b'} & \cb\ten_\ca\cb \ar[r] & 
\cb\ar[r] & 0}\ko
\end{equation}
where the map $b'$ takes $x\ten t\ten y$ to $xt\ten y -x\ten t y$ and the second
map is composition.
By taking the derived tensor product over $\cb^e$ with $\cb$, we find that
Hochschild homology of $\cb$ is computed by the cone over the induced morphism
\[
\xymatrix{ (\cb\ten_\ca\theta\ten_\ca\cb)\lten_{\cb^e} \cb \ar[r] & 
(\cb\ten_\ca\cb)\lten_{\cb^e}\cb.}
\]
Notice that $\theta$ is cofibrant over $\ca^e$ so that
$\cb\ten_\ca\theta\ten_\ca\cb=\theta\ten_{\ca^e}\cb^e$ is cofibrant over
$\cb^e$ and for the first term, we have
\[
 (\cb\ten_\ca\theta\ten_\ca\cb)\lten_{\cb^e} \cb =(\theta\ten_{\ca^e}\cb^e)\lten_{\cb^e}\cb
 =\theta\ten_{\ca^e}\cb.
\]
As for the second term, notice that $\cb$ is cofibrant over $\ca$ so that
we have 
\[
\cb\ten_\ca\cb=\cb\lten_\ca\cb=\ca\lten_{\ca^e}\cb^e
\] 
and therefore
\[
(\cb\ten_\ca\cb)\lten_{\cb^e}\cb = (\ca\lten_{\ca^e}\cb^e)\lten_{\cb^e}\cb = 
\ca\lten_{\ca^e}\cb.
\]
Now notice that in $\cb=\ca\oplus\theta\oplus(\theta\ten_\ca\theta)\oplus \cdots$,
all the summands are cofibrant over $\ca^e$ except the first one. Since we
are interested in {\em reduced} Hochschild homology, the first term does
not matter: Let us put $\ol{\cb}=\cb/\ca$. Then $\ol{\cb}$ is cofibrant over
$\ca^e$ and we have $\ca\lten_{\ca^e}\ol{\cb} = \ca\ten_{\ca^e}\ol{\cb}$.
So the reduced Hochschild homology of $\cb$ is computed by the cone over
the morphism
\[
\xymatrix{\theta\ten_{\ca^e} \cb \ar[r]^-b & \ca\ten_{\ca^e}\ol{\cb}.}
\]
It is not hard to check that the map $b$ is given by 
\[
t\ten u \mapsto 1_y \ten tu - (-1)^{|t||u|} 1_x \ten ut \ko
\]
where $t\in \theta(x,y)$ and $u\in\cb(y,x)$. Notice that its kernel has the summand
$\theta \ten_{\ca^e}\ca$. The construction shows that the cone over $b$
is the analogue of the complex $C^{\mathrm{small}}(T(V))=(V\ten A \to A)$ of section~3.1.1 of \cite{Loday98}, where $A=T(V)$. Now proceding further along this line, it is not
hard to check that the analogue of the complex
\[
\xymatrix{\ldots \ar[r]^-\gamma & V\ten A \ar[r]^-b & A \ar[r]^-\gamma & V\ten A \ar[r]^-b &
A \ar[r] & 0}
\]
computing the cyclic homology of $A=T(V)$ in Prop.~3.1.5 of \cite{Loday98}
is the sum total dg module of the periodic complex
\[
\xymatrix{\ldots \ar[r]^-\gamma & \theta \ten_{\ca^e}\cb \ar[r]^-b &
\ca\ten_{\ca^e}\ol{\cb} \ar[r]^-\gamma & \theta\ten_{\ca^e}\cb \ar[r]^-b &
\ca\ten_{\ca^e}\ol{\cb} \ar[r] & 0}\ko
\]
which computes the reduced cyclic homology of $\cb$ and where 
$\gamma$ is given by
\[
a \ten (t_1 \ldots t_n) \mapsto 
\sum_{i=1}^n \pm t_i \ten (t_{i+1} \ldots t_n a t_1 \ldots t_{i-1})\ko
\]
where the sign is given by the Koszul sign rule.
Let us now exhibit a canonical element in the reduced $(n-1)$th cyclic homology
of $\cb$. We have canonical quasi-isomorphisms
\[
\xymatrix{\theta \ten_{\ca^e} \ca & \theta\ten_{\ca^e} \bp\ca \ar[l] \ar[r] & \Sigma^{n-1}\Hom_{\ca^e}(\bp \ca,\bp\ca)} \ko
\]
where $\bp\ca\to\ca$ is an $\ca^e$-cofibrant resolution of $\ca$.
The identity on the right hand side corresponds to a Casimir element $c$ on the
left hand side. This element yields a canonical class of homological
degree $n-1$ in the summand
$\ca\ten_{\ca^e}\theta\cong\theta\ten_{\ca^e}\ca$ of the last term of the complex\[
\xymatrix{\ldots \ar[r]^-\gamma & \theta \ten_{\ca^e}\cb \ar[r]^-b &
\ca\ten_{\ca^e}\ol{\cb} \ar[r]^-\gamma & \theta\ten_{\ca^e}\cb \ar[r]^-b &
\ca\ten_{\ca^e}\ol{\cb} \ar[r] & 0}
\]
and thus a canonical class in $HC_{n-1}^{\mathrm{red}}(\cb)$. Under Connes'
map $HC_{n-1}^{\mathrm{red}}(\cb) \to HH_n^{\mathrm{red}}(\cb)$, which
is induced by $\gamma$, this class clearly corresponds to the class of the
Casimir element $c$ in the summand $\theta \ten_{\ca^e} \ca$ of 
$\theta\ten_{\ca^e}\cb$. It remains to be checked that this class
is non degenerate,\ie that it is taken to an isomorphism by the canonical maps
\[
\cb\lten_{\cb^e}\cb \iso \cb\lten_{\cb^e}\Theta^\vee \iso
\RHom_{\cb^e}(\Theta, \cb) \ko
\]
where $\Theta$ is the inverse dualizing complex $\cb^\vee$.
By the exact sequence (\ref{eq:small-resolution}), the $\cb$-bimodule $\cb$
is quasi-isomorphic to the cone over the map
\[
\xymatrix{\cb\ten_\ca\theta\ten_\ca\cb \ar[r]^-{b'} & \cb\ten_\ca\cb.}
\]
As we have already noted, the left hand term is $\cb^e$-cofibrant but the right
hand term is not. We choose a surjective $\ca^e$-cofibrant resolution $\bp\ca\to\ca$.
It yields a surjective quasi-isomorphism
\[
\xymatrix{\cb\ten_\ca\bp\ca\ten_\ca\cb \ar[r] & \cb\ten_\ca\cb}\ko
\]
where now the left hand term is $\cb^e$-cofibrant. We choose a lift
$\tilde{b}'$ of $b'$ along this quasi-isomorphism. We find that $\cb$ has as
$\cb^e$-cofibrant resolution the cone over the map
\[
\xymatrix{\cb\ten_\ca\theta\ten_\ca\cb \ar[r]^-{\tilde{b}'} &
\cb\ten_\ca\bp\ca\ten_\ca\cb}\ko
\]
which we can rewrite as
\[
\xymatrix{\theta\ten_{\ca^e}\cb^e \ar[r] & \bp\ca\ten_{\ca^e}\cb^e}.
\]
By applying $\Hom_{\cb^e}(?,\cb^e)$ to this cone and using the adjunction, 
we find that $\Theta$ is the cylinder over
\[
\xymatrix{\Hom_{\ca^e}(\bp\ca,\cb^e) \ar[r] & \Hom_{\ca^e}(\theta,\cb^e).}
\]
We know that the class of $c$ in $\theta\ten_{\ca^e}\ca\subset \theta\ten_{\ca^e}\cb$ yields a morphism from this cylinder shifted by $n$ degrees to the 
$\cb^e$-cofibrant resolution of $\cb$, \ie the cone over
\[
\xymatrix{\cb^e\ten_{\ca^e}\theta \ar[r] & \cb^e\ten_{\ca^e} \bp \ca.}
\]
This morphism is invertible in $\cd(\cb^e)$.
Indeed, it is not hard to see that in the left respectively right hand term, 
the class $c$ induces the canonical isomorphisms of $\cd(\cb^e)$
\[
\Hom_{\ca^e}(\bp\ca,\cb^e) \iso\cb^e\ten_{\ca^e}\Sigma^{1-n}\theta 
\quad\mbox{resp.}\quad
\Hom_{\ca^e}(\theta,\cb^e) \iso \cb^e\ten_{\ca^e}\Sigma^{1-n}\bp\ca.
\]
\end{proof}

\section{Deformed Calabi--Yau completions}
We keep the notations and assumptions of section~\ref{s:CY-completions}.
In particular, $\ca$ is assumed to be homologically smooth.
Let $c$ be a class in the Hochschild homology $HH_{n-2}(\ca)$.
Via the canonical isomorphism
\[
HH_{n-2}(\ca) \iso \Hom_{\cd(\ca^e)}(\theta,\Sigma\ca)\ko
\]
the class $c$ may be lifted to a closed degree~$1$ bimodule morphism
$\delta: \theta\to\ca$. The {\em deformed Calabi--Yau completion $\Pi_n(\ca,c)$}
was defined in section~5 of \cite{Keller11b} as obtained from
$\Pi_n(\ca)=T_\ca(\theta)$ by replacing the differential $d$ with the unique
derivation of the tensor category extending 
\[
\left[\begin{array}{cc} d & 0 \\ \delta & d \end{array} \right] : 
\theta\oplus \ca \to \theta\oplus\ca.
\]
Up to isomorphism, it only depends on $\ca$ and the class $c$. It was claimed in
Theorem~5.2 of \cite{Keller11b}, that the deformed Calabi--Yau 
completion is always $n$-Calabi--Yau. A counter-example to this
claim is given by W.~K.~Yeung at the end of section~3.3 in \cite{Yeung16}.
As pointed out by Yeung, a sufficient condition for the deformed
$n$-Calabi--Yau completion to be $n$-Calabi--Yau is suggested by
the work of Thanhoffer de V\"olcsey--Van den Bergh \cite{ThanhofferVandenBergh12}:
It suffices that the class $c$ lifts to the negative cyclic homology $HN_{n-2}(\ca)$.
The following is Theorem~3.17 of \cite{Yeung16}. 

\begin{theorem}[Yeung] Suppose that $\ca$ is finitely cellular and $c$ is a class in
$HH_{n-2}(\ca)$. Any lift $\tilde{c}$ of $c$ to the negative cyclic homology 
$HN_{n-2}(\ca)$ determines an $n$-Calabi--Yau structure on the
deformed Calabi--Yau completion $\Pi_n(\ca,c)$.
\end{theorem}

Notice that the assumptions do hold for Ginzburg dg categories.
One would expect that, like theorem~\ref{thm:cy-completion}, this theorem
should generalize to arbitrary homologically smooth dg categories.

\section{Correction of some inaccuracies in ``Deformed Calabi--Yau completions"}
\label{s:inaccuracies}

In section~3.6, one should point out in addition
that the assumption that $Q$ is projective over $k$ implies that it is 
projective as an $\cR$-bimodule (since $\cR$ is discrete, so
is $\cR^e$), which in turn implies that
the filtration (3.6.1) is split in the category of $\cR$-bimodules.

At the end of the proof of Theorem~4.8, it is claimed that the transpose
conjugate of $\tilde{\lambda}$ maps to $\rho$. This is not true. However,
the image of $\tilde{\lambda}$ and $\rho$ define the same class in homology.

The statement of Theorem~6.1 is incorrect. One has to assume in addition
that the set of `minimal relations' $R$ generates the
ideal $I$ (from the definition of $R$, it only follows that 
$R$ topologically generates the $J$-adic completion of $I$).
Moreover, the condition 3) in the proof of the theorem should be
replaced with
\begin{quote}
for all $n\geq 1$, the differential $d$ maps $V^{-n-1}$ to $T_n$ 
and induces an isomorphism from $V^{-n-1}$ onto the head
of the $H^0(T_n)$-bimodule $H^{-n}(T_n)$, where $T_n$
denotes the dg category $T_{\mathcal{R}}(V^0\oplus\cdots \oplus V^{-n})$.
\end{quote}

\section*{Acknowledgments}

I am very grateful to W.~K.~Yeung for providing a counter-example to Theorem~5.2
of \cite{Keller11b} and indicating a way to correct the error. I thank
Xiaofa Chen, Martin Kalck, Kai Wang and Dong Yang for pointing out
the inaccuracies corrected in section~\ref{s:inaccuracies}.


\begin{thebibliography}{10}

\bibitem{CrawleyBoevey00}
William Crawley-Boevey, \emph{On the exceptional fibres of {K}leinian
  singularities}, Amer. J. Math. \textbf{122} (2000), no.~5, 1027--1037.

\bibitem{Ginzburg06}
Victor Ginzburg, \emph{{Calabi-Yau} algebras}, arXiv:math/0612139v3 [math.AG].

\bibitem{IkedaQiu18}
Akishi Ikeda and Yu~Qiu, \emph{$\mathbb{X}$-stability conditions on
  {C}alabi--{Y}au-$\mathbb{X}$ categories and twisted periods},
  arXiv:1807.00469 [math.AG].

\bibitem{Keller11b}
Bernhard Keller, \emph{Deformed {C}alabi--{Y}au completions}, Journal f{\"u}r
  die reine und angewandte Mathematik (Crelles Journal) \textbf{654} (2011),
  125--180, with an appendix by Michel~Van den Bergh.

\bibitem{KontsevichSoibelman06}
Maxim Kontsevich and Yan Soibelman, \emph{{Notes on A-infinity algebras,
  A-infinity categories and non-commutative geometry. I}},
  arXiv:math.RA/0606241.

\bibitem{KontsevichVlassopoulos13}
Maxim Kontsevich and Yannis Vlassopoulos, \emph{Weak {C}alabi-{Y}au algebras},
  talk at a Conference on Homological Mirror Symmetry, University of Miami,
  2013, available at \verb"https://math.berkeley.edu/~auroux/miami2013.html".

\bibitem{Loday98}
Jean-Louis Loday, \emph{Cyclic homology}, second ed., Grundlehren der
  Mathematischen Wissenschaften [Fundamental Principles of Mathematical
  Sciences], vol. 301, Springer-Verlag, Berlin, 1998, Appendix E by Mar\'\i a
  O. Ronco, Chapter 13 by the author in collaboration with Teimuraz Pirashvili.

\bibitem{VandenBergh15}
Michel Van~den Bergh, \emph{Calabi-{Y}au algebras and superpotentials}, Selecta
  Math. (N.S.) \textbf{21} (2015), no.~2, 555--603.

\bibitem{ThanhofferVandenBergh12}
Michel Van~den Bergh and Louis de~Thanhoffer~de V\"olcsey, \emph{Calabi--{Y}au
  deformations and negative cyclic homology}, arXiv:1201.1520 [math.RA].

\bibitem{Yeung16}
Wai-kit Yeung, \emph{Relative {C}alabi--{Y}au completions}, arXiv:1612.06352
  [math.RT].

\bibitem{Yeung18}
\bysame, \emph{Weak {C}alabi--{Y}au structures and moduli of representations},
  arXiv:1802.05398 [math.AG].

\end{thebibliography}

\def\cprime{$'$} \def\cprime{$'$}
\providecommand{\bysame}{\leavevmode\hbox to3em{\hrulefill}\thinspace}
\providecommand{\MR}{\relax\ifhmode\unskip\space\fi MR }
\providecommand{\MRhref}[2]{%
  \href{http://www.ams.org/mathscinet-getitem?mr=#1}{#2}
}
\providecommand{\href}[2]{#2}

\end{document}